\documentclass[reqno]{amsart}

\usepackage{amsfonts}
\usepackage{amsmath,amssymb}

 \DeclareMathOperator{\E}{{\mathbb{E}}}
 \DeclareMathOperator{\Prob}{{\mathbb{P}}}

\newcommand{\ii}{{\mathrm{i}}}
\newcommand{\dd}{{\mathrm{d}}}

\newtheorem{thm}{Theorem}
\newtheorem{cor}[thm]{Corollary}
\newtheorem{lem}[thm]{Lemma}
\newtheorem{proposition}[thm]{Proposition}
\theoremstyle{remark}
\newtheorem*{rem}{Remark}
\theoremstyle{definition}
\newtheorem*{example}{Example}

\begin{document}
\title[Zeros of random polynomials]{The zeros of random polynomials cluster
uniformly near the unit circle}

\author{C.P. Hughes}
\address{Department of Mathematics, University of York,
York, YO10 5DD, U.K.}
 \email{ch540@york.ac.uk}

\author{A. Nikeghbali}
 \address{Institut f\"ur Mathematik,
Universit\"at Z\"urich, Winterthurerstrasse 190, CH-8057 Z\"urich,
Switzerland}
 \email{ashkan.nikeghbali@math.unizh.ch}

\date{3 June 2007}

\begin{abstract}
In this paper we deduce a universal result about the asymptotic
distribution of roots of random polynomials, which can be seen as
a complement to an old and famous result of Erd\H{o}s and Turan.
More precisely, given a sequence of random polynomials, we show
that, under some very general conditions, the roots tend to
cluster near the unit circle, and their angles are uniformly
distributed. The method we use is deterministic: in particular, we
do not assume independence or equidistribution of the coefficients
of the polynomial.
\end{abstract}

\maketitle

\section{Introduction}\label{sect:intro}

In this paper, we are interested in the uniform concentration near
the unit circle of roots of polynomials.

Let $\left(P_{N}(Z)\right)_{N\geq 1}$ be a sequence of
polynomials. Denote the zeros of $P_N(Z)$ by $z_1,\dots,z_N$. Let
\begin{equation}\label{eq:defn_nu_zer}
\nu _{N}\left( \rho \right) :=\#\left\{ z_{k}\ :\ 1-\rho \leq
\left\vert z_{k}\right\vert \leq \frac{1}{1-\rho }\right\}
\end{equation}
be the number of zeros of $P_{N}(Z) $ lying in the annulus bounded
by $1-\rho $ and $\frac{1}{1-\rho }$, where $0\leq \rho \leq 1$,
and let
\begin{equation}\label{eq:defn_nu_ang}
\nu _{N}\left( \theta ,\phi \right) :=\#\left\{ z_{k}\ :\ \theta
\leq \arg \left( z_{k}\right) <\phi \right\}
\end{equation}
be the number of zeros of $P_{N}(Z) $ whose argument lies between
$\theta $ and $\phi $, where $0\leq \theta <\phi \leq 2\pi $.

We shall say that the zeros cluster uniformly around the unit
circle if for all fixed $0<\rho<1$,
\begin{equation}\label{propzer}
\lim_{N\rightarrow \infty }\frac{1}{N}\nu _{N}\left( \rho \right)
=1
\end{equation}
and
\begin{equation}\label{propang}
\lim_{N\rightarrow \infty }\frac{1}{N}\nu _{N}\left( \theta ,\phi
\right) = \frac{\phi -\theta }{2\pi }
\end{equation}

The purpose of this paper is to find a general but simple
condition for when the zeros cluster uniformly around the unit
circle.

\begin{thm}\label{thm:thm1}
Let $\left(P_{N}(Z)\right)$ be a sequence of polynomials, with
\begin{equation*}
P_{N}(Z) =\sum_{k=0}^{N}a_{N,k}Z^{k},
\end{equation*}
such that $ a_{N,0}a_{N,N}\neq0 $ for all $N$. Let
\begin{equation*}
L_{N}\left(P_{N}\right)=\log \left( \sum_{k=0}^{N}|a_{N,k}|\right)
-\frac{1}{2}\log |a_{N,0}|- \frac{1}{2}\log |a_{N,N}| .
\end{equation*}
If $$L_{N}\left(P_{N}\right)=o\left(N\right),$$ then the zeros of
this sequence cluster uniformly near the unit circle, i.e. for all
$0<\rho<1$,
\begin{equation*}
\lim_{N\rightarrow \infty }\frac{1}{N}\nu _{N}\left(\rho\right) =1,
\end{equation*} and for all $0\leq
\theta<\phi\leq 2\pi$,
\begin{equation*}
\lim_{N\rightarrow \infty } \frac{1}{N}\nu_{N}\left( \theta ,\phi
\right)  =\frac{\phi -\theta }{2\pi }.
\end{equation*}
where $\nu _{N}\left(\rho\right)$ and $\nu_{N}\left( \theta ,\phi
\right)$ are defined in \eqref{eq:defn_nu_zer} and
\eqref{eq:defn_nu_ang} respectively.
\end{thm}

The second part of our theorem, on $\nu_N(\theta,\phi)$, follows
from the celebrated result of Erd\H{o}s and Turan
\cite{Erdosturan} on the distribution of roots of polynomials:

\begin{thm}[Erd\H{o}s-Turan]
\label{ErdosTuran} Let $\left( a_{k}\right)_{0\leq k \leq N}$ be a
sequence of complex numbers such that $a_{0}a_{N}\neq 0$, and let
\begin{equation*}
P(Z) = \sum_{k=0}^N a_k Z^k \ .
\end{equation*}
Define
\begin{equation}\label{LN}
L_{N}\left(P\right)=\log \sum_{k=0}^{N}\left\vert a_{k}\right\vert
-\frac{1}{2}\log \left\vert a_{0}\right\vert -\frac{1}{2} \log
\left\vert a_{N}\right\vert.
\end{equation}
Then
\begin{equation*}
\left\vert \frac{1}{N}\nu _{N}\left(\theta,\phi\right) -\frac{\phi
-\theta }{2\pi }\right\vert ^{2}\leq
\frac{C}{N}L_{N}\left(P\right)
\end{equation*}
for some constant $C$, where $\nu_N(\theta,\phi)$ is defined in
\eqref{eq:defn_nu_ang}.
\end{thm}
The above theorem shows that if $L_{N}\left(P\right)$ is small
compared to the degree $N$, then the angles of the roots are nearly
uniformly distributed, and that is precisely the reason why this
theorem has been extensively used to prove asymptotic uniform
concentration near the unit circle of the roots of some families of
random polynomials.

In this paper, we prove a natural complement to this result.
\begin{thm}\label{Hughes-NI}
Let $\left( a_{k}\right)_{0\leq k \leq N}$ be a sequence of
complex numbers such that $a_{0}a_{N}\neq 0$, and let
\[
P(Z) = \sum_{k=0}^N a_k Z^k.
\]
Let $L_N(P)$ be defined as in \eqref{LN}. Then for $0<\rho<1$,
\begin{equation*}
\left(1 - \frac{1}{N}\nu_N(\rho)\right) \leq \frac{2}{N\rho}
L_N(P)
\end{equation*}
where $\nu_N(\rho)$ is defined in \eqref{eq:defn_nu_zer}.
\end{thm}

These two theorems on deterministic polynomials give a sufficient
condition for the roots of random polynomials to cluster uniformly
around the unit circle, and we will show how results of \v{S}paro
and \v{S}ur \cite{SparoSur}, Arnold \cite{Arnold}, and Shmerling
and Hochberg \cite{Schmer} follow as corollaries of this. Indeed,
we show that some of their conditions on the coefficients of the
polynomials can be dropped.

More precisely, let $\left( \Omega ,\mathcal{F},\Prob\right) $ be
a probability
space, on which a random array $(a_{N,k})_{\substack{ N\geq 1  \\
0\leq k\leq N}}$ is defined. Consider now the sequence of random
polynomials $\left(P_{N}(Z)\right)$, with
\begin{equation}\label{tableau}
P_{N}(Z) =\sum_{k=0}^{N}a_{N,k}Z^{k}.
\end{equation}

The asymptotics for roots of random polynomials in the complex plane
have already been studied in the special case where
\begin{equation*}
P_{N}(Z) =\sum_{k=0}^{N}a_{k}Z^{k};
\end{equation*}
which corresponds to the special array $\left( a_{N,k}\right)
_{k\leq N}=\left( a_{0},a_{1},\ldots ,a_{N}\right) $, and it is
not our intention to give a full historical account here (see
\cite{Bharucha}, \cite{Farahmand}, \cite{EdelKost} for more
details and references), but we shall rather mention the papers of
\v{S}paro and \v{S}ur \cite{SparoSur}, Arnold \cite{Arnold}, and
Shmerling and Hochberg \cite{Schmer} which contain the most
general results about uniform clustering near the unit circle.
\v{S}paro and \v{S}ur \cite{SparoSur} have considered i.i.d.
complex coefficients $\left( a_{n}\right) _{n\geq 0}$ and have
shown that under some integrability conditions, the zeros of such
sequences cluster uniformly near the unit circle, with convergence
in \eqref{propzer} and \eqref{propang} holding in probability.
Arnold \cite{Arnold} improved this result and proved that the
convergence holds almost surely and in the $p$\textsuperscript{th}
mean if the moduli of $a_{k}$ are equidistributed (plus some
integrability conditions). Recently, Shmerling and Hochberg
\cite{Schmer} have obtained stronger results: they have shown that
the condition on equidistribution can be dropped if $\left(
a_{n}\right) _{n\geq 0}$ is a sequence of independent variables
which have continuous densities $f_{n}$ which are uniformly
bounded in some neighbourhood of the origin with finite means $\mu
_{n}$ and standard deviations $\sigma _{n}$\ that satisfy the
condition
\begin{gather*}
\max \left\{ \limsup_{n\rightarrow \infty }\sqrt[n]{|\mu _{n}|}\ ,\
\limsup_{n\rightarrow \infty }\sqrt[n]{|\sigma _{n}|}\right\} =1, \\
\Prob\left\{ a_{0}=0\right\} =0 .
\end{gather*}

The authors above mentioned consider separately the cases of
\eqref{propzer} and \eqref{propang}. They prove \eqref{propang}
using Theorem \ref{ErdosTuran}, and prove \eqref{propzer} using
techniques from random power series. In particular, to prove
\eqref{propang}, they are able to show that
$\frac{L_{N}\left(P\right)}{N}\rightarrow 0$, as
$N\rightarrow\infty$ and this is only a few lines  proof, while
the techniques used to prove \eqref{propzer} are more
sophisticated. Thus a side benefit of Theorem~\ref{Hughes-NI} is
that once one proves the uniform distribution of the angles using
Theorem~\ref{ErdosTuran} of Erd\H{o}s and Turan, then one actually
proves uniform clustering near the unit circle considerably
simplifying the arguments in the proofs on uniform clustering in
\cite{Arnold} and \cite{Schmer}. Moreover, as we shall see, our
results also give us an estimate for the rate of clustering (that
is, how quickly $\rho\to0$ as a function of $N$).

The historical results we have mentioned above do not apply to the
more general case of sequences of random polynomials of the form
\eqref{tableau} we are dealing with, and they also do not lead to
anything interesting in the case of deterministic sequences of
polynomials (the asymptotic study of roots sequences of
deterministic polynomials occurs in problems of equidistribution
of algebraic integers, for example \cite{Bilu}). Moreover, these
results do not cover the cases where the coefficients are
dependent with different distributions: for example in the
semiclassical approximations for multidimensional quantum systems,
one needs to locate roots of high degree random polynomials with
dependent and non identically distributed coefficients \cite{BBL}
(we should point out that the authors in \cite{BBL} have observed
the uniform clustering in the special case of self-reciprocal
polynomials with complex Gaussian coefficients with finite
variance).

The aim of this paper is to show that the phenomenon of uniform
concentration of zeros around the unit circle is universal, in the
sense that no independence or equidistribution on the coefficients
is required, but only conditions on their size. Our method, based
on elementary complex analysis, reduces both convergences
\eqref{propzer} and \eqref{propang} to the same problem, namely
showing that $L_N(P)$, defined in \eqref{LN}, is small compared to
the degree $N$ of the polynomial $P$, thus complementing
Theorem~\ref{ErdosTuran} of Erd\H{o}s and Turan.

The structure of this paper is as follows: In
section~\ref{sect:estimates} we prove our main theorems  and then in
section~\ref{sect:clustering} we  use them to deduce clustering of
zeros for general sequences of random polynomials.

\section{Basic estimates}

\label{sect:estimates} For $N\geq 1$, let $\left(
a_{k}\right)_{0\leq k\leq N}$ be a sequence of complex numbers
satisfying $a_{0} a_{N} \neq 0$. {From} this sequence construct
the polynomial
\begin{equation*}
P_{N}\left(Z\right) =\sum_{k=0}^{N} a_{k}Z^{k},
\end{equation*}
and denote its zeros by $z_{i}$ (where $i$ ranges from 1 to $N$). For $0\leq
\rho\leq 1$, we are interested in estimates for
\begin{eqnarray*}
\widetilde{\nu}_{N}(\rho) &=& \#\left\{ z_{j} \ , \ |z_{j}| <
1-\rho\right\} \\
\overline{\nu}_{N}(\rho) &=& \#\left\{ z_{j} \ , \ |z_{j}| >
\frac{1}{1-\rho} \right\} \\
\nu_{N}(\rho) &=& \#\left\{ z_j \ , \ 1-\rho \leq |z_j| \leq \frac{1}{1-\rho}
\right\}
\end{eqnarray*}
which counts the number of zeros of the polynomial $P_{N}(Z) $
which lie respectively inside the open disc of radius $1-\rho$,
outside the closed disc of radius $1/(1-\rho)$, and inside the
closed annulus bounded by circles of radius $1-\rho$ and
$1/(1-\rho)$.

\begin{thm}
\label{analytictheorem} For $N\geq 1$, let $(a_{k})_{0\leq k \leq
N}$ be a sequence of complex numbers which satisfy $a_{0} a_{N}
\neq 0$. Then, for $0<\rho<1$
\begin{equation}  \label{essentialineq1}
\frac{1}{N} \widetilde{\nu}_{N}\left( \rho\right) \leq
\frac{1}{N\rho} \left( \log \left(\sum_{k=0}^N |a_{k}|\right) -
\log |a_{0}| \right) ,
\end{equation}
\begin{equation}  \label{essentialineq2}
\frac{1}{N} \overline{\nu}_{N}(\rho) \leq \frac{1}{N\rho} \left(
\log \left(\sum_{k=0}^n |a_{k}|\right) - \log |a_{N}|\right)
\end{equation}
and
\begin{equation}  \label{essentialineq3b}
\left(1 - \frac{1}{N}\nu_N(\rho)\right) \leq \frac{2}{N\rho} \left( \log
\left( \sum_{k=0}^N |a_{k}| \right) - \frac{1}{2} \log |a_{0}| - \frac{1}{2}
\log |a_{N}| \right) .
\end{equation}
\end{thm}

\begin{proof}
An application of Jensen's formula (see, for example, \cite{Lang}
page 341) yields
\begin{equation*}
\frac{1}{2\pi }\int_{0}^{2\pi }\log \left\vert P_{N}(e^{\ii
\varphi })\right\vert \;\dd\varphi -\log
|P_{N}(0)|=\sum_{|z_{i}|<1}\log \frac{1}{|z_{i}|}
\end{equation*}
where the sum on the right hand side is on zeros lying inside the open unit
disk. We have the following minorization for this sum:
\begin{align*}
\sum_{|z_{i}|<1}\log \frac{1}{|z_{i}|}& \geq \sum_{|z_{i}|<1-\rho }\log
\frac{1}{|z_{i}|} \\
& \geq \rho \widetilde{\nu }_{N}\left( \rho \right)
\end{align*}
since if $0\leq \rho \leq 1$, then for all $|z_{i}|\leq 1-\rho $,
$\log (1/|z_{i}|)\geq \rho $, and by definition there are
$\widetilde{\nu }_{N}(\rho )$ such terms in the sum.

We also have the following trivial upper bound
\begin{equation*}
\max_{\varphi \in [0,2\pi]} |P_N(e^{\ii\varphi})| \leq
\sum_{k=0}^N |a_k| ,
\end{equation*}
and so
\begin{eqnarray*}
\rho \widetilde{\nu }_{N}\left( \rho \right) &\leq &\frac{1}{2\pi}
\int_{0}^{2\pi }\log |P_{N}( e^{\ii \varphi })| \; \dd \varphi -\log |a_{0}| \\
&\leq& \log \left( \sum_{k=0}^{N} |a_{k}| \right) -\log | a_{0}|
\end{eqnarray*}
which gives equation \eqref{essentialineq1}.

To estimate the number of zeros lying outside the closed disc of
radius $\left( 1-\rho\right)^{-1} $, note that if $z_{0}$ is a
zero of the polynomial $P_{N}(Z) =\sum_{k=0}^{N}a_{k}Z^{k}$, then
$1/z_{0}$ is a zero of the polynomial $Q_{N}(Z) :=Z^{N}P_{N}\left(
\frac{1}{Z} \right) =a_{N}+a_{N-1}Z+\ldots +a_{0}Z^{N}$.
Therefore, the number of zeros of $P_{N}(Z)$ outside the closed
disc of radius $1/(1-\rho)$ equals the number of zeros of
$Q_{N}(Z)$ inside the open disc of radius $1-\rho$. Therefore,
from \eqref{essentialineq1} we get
\begin{equation*}
\frac{1}{N} \overline{\nu}_{N}(\rho) \leq \frac{1}{N\rho} \left(
\log \left(\sum_{k=0}^N |a_{k}|\right) - \log |a_{N}| \right)
\end{equation*}
which gives equation~\eqref{essentialineq2}.

Since
\begin{equation*}
N - \nu_N(\rho) =\widetilde\nu_N(\rho) + \overline\nu_N(\rho)
\end{equation*}
we immediately get \eqref{essentialineq3b}.
\end{proof}

Theorem~\ref{Hughes-NI} is merely a restatement of equation
\eqref{essentialineq3b}. Theorem~\ref{thm:thm1} follows from
combining this with the result of Erd\H{o}s and Turan,
Theorem~\ref{ErdosTuran}.

\begin{rem}
Note that if $a_{k}\mapsto \lambda a_{k}$ for some $\lambda \neq 0$, then
the zeros of $P_{N}(Z)$ are unchanged, and
\begin{multline*}
\log \left( \sum_{k=0}^{N}|\lambda a_{k}|\right) -\frac{1}{2}\log |\lambda
a_{0}|-\frac{1}{2}\log |\lambda a_{N}| \\
=\log \left( \sum_{k=0}^{N}|a_{k}|\right) -\frac{1}{2}\log
|a_{0}|-\frac{1}{2}\log |a_{N}|
\end{multline*}
so, in some sense, this is a natural function to control the location of the
zeros.
\end{rem}

\section{Uniform clustering results for roots of random polynomials}

\label{sect:clustering}

We require no independence restriction on our random variables. We only
assume that
\begin{equation} \label{cond1}
\Prob\left\{ a_{N,0}=0\right\} =0  
\end{equation}
and
\begin{equation} \label{cond1bis}
\Prob\left\{ a_{N,N}=0\right\} =0,  
\end{equation}
for all $N$.

\subsection{The main theorem for random polynomials}

\begin{thm}
\label{generalclustering} For $N\geq 1$, let $(a_{N,k})_{0\leq
k\leq N}$ be an array of random complex numbers such that
$\Prob\left\{ a_{N,0}=0\right\} =0 $ and $\Prob\left\{
a_{N,N}=0\right\} =0$ for all $N$. Let
\begin{equation}\label{eq:defn F_N}
L_{N}=\log \left( \sum_{k=0}^{N}|a_{N,k}|\right) -\frac{1}{2}\log
|a_{N,0}|- \frac{1}{2}\log |a_{N,N}| .
\end{equation}

If
\begin{equation}\label{eq:assumption on F_N}
\E\left[ L_{N}\right] =o(N)\quad \text{ as }N\rightarrow \infty
\end{equation}
then there exists a positive function $\alpha _{N}$ satisfying
$\alpha _{N}=o(N)$ such that the zeros of the random polynomial
\begin{equation*}
P_{N}(Z)=\sum_{k=0}^{N}a_{N,k}Z^{k}
\end{equation*}
satisfy
\begin{equation*}
\lim_{N\rightarrow \infty }\E\left[ \frac{1}{N}\nu _{N}\left( \frac{\alpha
_{N}}{N}\right) \right] =1
\end{equation*}
and
\begin{equation*}
\lim_{N\rightarrow \infty }\E\left[ \frac{1}{N}\nu _{N}\left( \theta ,\phi
\right) \right] =\frac{\phi -\theta }{2\pi } .
\end{equation*}
In fact the convergence also holds in probability and in the $p$\textsuperscript{th} mean, for all positive $p$.

Furthermore, if there exists a (deterministic) positive function $\alpha
_{N} $ satisfying $\alpha _{N}\leq N$ for all $N$, such that
\begin{equation}
L_{N}=o(\alpha _{N})\quad \text{almost surely} \label{eq:almost sure
assumption on F_N}
\end{equation}
then both convergences hold almost surely (and also in the
$p$\textsuperscript{th} mean, for all positive $p$).
\end{thm}
\begin{proof}
The convergence in mean for $\nu _{N}(\alpha _{N}/N)$ is a
consequence of \eqref{essentialineq3b}. We have
\begin{equation*}
1-\E\left[ \frac{1}{N}\nu _{N}\left( \frac{\alpha _{N}}{N}\right)
\right] \leq \frac{2}{\alpha _{N}}\E\left[ L_{N}\right] .
\end{equation*}
Therefore we see that the result follows for any positive function
$\alpha _{N}$ satisfying $\alpha _{N}\leq N$ for all $N$ such that
$\E\left[ L_{N} \right] /\alpha _{N}\rightarrow 0$, and such a
function exists by assumption \eqref{eq:assumption on F_N}, for
example
\[
\alpha_N = N \min\left\{1,\sqrt{\frac{\E[L_N]}{N}}\right\} .
\]

Similarly from Theorem \ref{ErdosTuran} and \eqref{eq:assumption on
F_N} we have that
\begin{align*}
\E\left[ \left\vert \frac{1}{N}\nu _{N}\left( \theta ,\phi \right)
-\frac{\phi -\theta }{2\pi }\right\vert ^{2}\right] & \leq
\frac{C}{N}\E\left[ L_{N}\right] \\
& =o(1) .
\end{align*}

Note that the mean square convergence implies convergence in the mean, as in
the theorem, and also convergence in probability. Note further, that since
the random variables are uniformly bounded ($0\leq \frac{1}{N}\nu
_{N}\left(\theta,\phi\right) \leq 1$), mean convergence implies convergence
in the $p$\textsuperscript{th} mean for all positive $p$.

In the same way, the almost sure convergence of $\frac{1}{N}\nu
_{N}(\alpha _{N}/N)$ and $\frac{1}{N}\nu _{N}(\theta ,\phi )$
follows immediately from \eqref{essentialineq3b} and Theorem
\ref{ErdosTuran}, using \eqref{eq:almost sure assumption on F_N}.
\end{proof}
In the following subsections, we shall give some sufficient
conditions, which are easy to check, for Theorem
\ref{generalclustering} to hold. We first consider the case of
general sequences of random polynomials for which there exist no
previous results to our knowledge; then we deal with the classical
random polynomials.

\subsection{General sequences of random polynomials}

\begin{proposition}
\label{resultconvprincip} Let $(a_{N,k})$ be an array of random
complex numbers which satisfy \eqref{cond1} and \eqref{cond1bis}.
Assume that $\E \left[ \log |a_{N,0}|\right] =o(N)$, and $\E\left[
\log |a_{N,N}|\right] =o(N)$, and that there exists a fixed
$0<s\leq 1$ such that
\begin{equation*}
\sum_{k=0}^{N}\E\left[ |a_{N,k}|^{s}\right] = \exp(o(N)) .
\end{equation*}
Then there exists a positive sequence $(\alpha_N)$ subject to
$\alpha _{N}=o(N)$ such that
\begin{equation*}
\lim_{N\rightarrow \infty }\E\left[ \frac{1}{N}\nu _{N}\left( \frac{\alpha
_{N}}{N}\right) \right] =1
\end{equation*}
and
\begin{equation*}
\lim_{N\rightarrow \infty }\E\left[ \left\vert \frac{1}{N}\nu _{N}\left(
\theta ,\phi \right) -\frac{\phi -\theta }{2\pi }\right\vert \right] =0 .
\end{equation*}
\end{proposition}

\begin{proof}
Since $0<s\leq 1$ we have the following concavity inequality:
\begin{equation*}
\E\left[ \log \left( \sum_{k=0}^{N}|a_{N,k}|\right) \right] \leq
\frac{1}{s}\log \left( \sum_{k=0}^{N}\E\left[ |a_{N,k}|^{s}\right]
\right) .
\end{equation*}
By assumption, the RHS is $o(N)$. Therefore $L_{N}$, defined in
\eqref{eq:defn F_N}, satisfies $L_{N}=o(N)$, and the result
follows from Theorem~\ref{generalclustering}.
\end{proof}

\begin{rem}
The Proposition shows that under some very general conditions
(just some conditions on the size of the expected values of the
modulus of the coefficients), without assuming any independence or
equidistribution condition, the zeros of random polynomials tend
to cluster uniformly near the unit circle. We can also remark that
we do not assume that our coefficients must have density
functions: they can be discrete-valued random variables.
\end{rem}

Now we give two examples which could not be dealt with the previous
results available in the literature.
\begin{example}
Let $a_{N,k}$ be random variables distributed according to the Cauchy
distribution with parameter $N(k+1)$. The first moment does not exist but
some fractional moments do, and in particular we have for $0\leq s<1$
\begin{align*}
\E\left[ |a_{N,k}|^{s}\right] & =\frac{N(k+1)}{\pi }\int_{-\infty
}^{\infty }\frac{|x|^{s}}{x^{2}+N^{2}(k+1)^{2}}\;\dd x \\
& =\frac{1}{\pi }N^{s}(k+1)^{s}\Gamma
(\frac{1}{2}+\frac{s}{2})\Gamma (\frac{1}{2}-\frac{s}{2}) .
\end{align*}
Moreover,
\begin{equation*}
\E\left[ \log |a_{N,k}|\right] =\log (N(k+1)) .
\end{equation*}
Hence we can apply Proposition ~\ref{resultconvprincip} and deduce that the
zeros of the sequence of random polynomials with coefficients $(a_{N,k})
_{\substack{ N\geq 1  \\ 0\leq k\leq N}}$ where $a_{N,k}$ are chosen from
the Cauchy distribution with parameter $N(k+1)$ cluster uniformly around the
unit circle.
\end{example}

\begin{example}
For each $N$, let $\left( a_{N,k}\right) $, $0\leq k\leq N$, be
discrete random variables taking values in $\left\{ \pm 1,\ldots
,\pm N\right\} $, not necessarily having the same distribution.
Then for any positive function $\alpha_N \leq N$ such that
$\alpha_N / \log N \to \infty$ we have
\begin{eqnarray*}
\lim_{N\rightarrow \infty }\frac{1}{N}\nu_{N}\left( \frac{\alpha_N}{N} \right)
&=&1,\quad a.s. \\
\frac{1}{N}\nu _{N}\left( \theta ,\phi \right) &\rightarrow &\frac{\phi
-\theta }{2\pi },\quad a.s.
\end{eqnarray*}
This follows from Theorem~\ref{generalclustering}, since for any
choice of the $a_{N,k}$ we have the deterministic bound
\begin{align*}
L_N &= \log\left(\sum_{k=0}^N |a_{N,k}| \right) -
\frac{1}{2}\log|a_{N,0}|
-\frac{1}{2}\log|a_{N,N}| \\
&\leq \log((N+1) N) < 2\log(N+1)
\end{align*}
with this choice of $\alpha_N$, $L_N = o(\alpha_N)$ a.s.

As a special case, we have the well known random polynomials
$\sum_{k=0}^{N}\mu _{k}Z^{k}$, where $\mu _{k}=\pm 1$, with
probabilities $p$ and $\left( 1-p\right) $. Moreover, we have from
the Markov inequality, the following rate for the convergence in
probability: There exists a constant $C$ such that
\begin{eqnarray*}
\Prob\left\{ \left( 1-\frac{1}{N}\nu _{N}\left( \frac{\alpha
_{N}}{N}\right) \right) >\varepsilon \right\} &\leq
&\frac{1}{\varepsilon }\frac{C\log N}{\alpha _{N}} \\
\Prob\left\{ \left\vert \frac{1}{N}\nu _{N}\left( \theta ,\phi
\right) -\frac{\phi -\theta }{2\pi }\right\vert >\varepsilon
\right\} &\leq &\frac{1}{\varepsilon ^{2}}\frac{C\log N}{N}
\end{eqnarray*}
for any fixed $\varepsilon >0$.
\end{example}

\subsection{Classical Random Polynomials}

\label{sect:classical random polynomials}

Let us now consider the special, but important, case of the classical random
polynomials as mentioned in the first section, that is
\begin{equation} \label{poly_simple}
P_{N}(Z) =\sum_{k=0}^{N}a_{k}Z^{k} . 
\end{equation}

These polynomials have been extensively studied (see, for example,
\cite{Bharucha} or \cite{Farahmand} for a complete account).

The results of the previous section take a simpler form in the
special case of random polynomials of the form
\eqref{poly_simple}. The conditions \eqref{cond1} and
\eqref{cond1bis} become
\begin{equation*}
\Prob\left\{ a_{N}=0\right\} =0,\text{ for all } N\geq 0 .
\end{equation*}

In this more special case, we can deal more easily with almost
sure convergence, which in our framework is the strongest
convergence. We will proceed to give some simple-to-check
sufficient conditions on the moments of $a_k$ to ensure
Theorem~\ref{generalclustering} holds.

\begin{thm}\label{thm:random polys}
Let $\left( a_{k}\right) _{k\geq 0}$ be a sequence of complex
random variables. Assume that
\begin{itemize}
\item There exists some $s>0$ such that if $\lambda_k := \E\left[
|a_{k}|^{s} \right] $ then $\lambda_k < \infty$ for all $k$, and
\begin{equation*}
\limsup_{k\to\infty}\left( \lambda _{k}\right) ^{1/k}  = 1 .
\end{equation*}
\item There exists some $t>0$ such that for all $k$,
\begin{equation}\label{eq:condition on neg power}
\xi_{k} := \E\left[ \frac{1}{|a_{k}|^{t}} \right] <\infty
\end{equation}
and
\begin{equation*}
\limsup_{k\to\infty} \left(\xi_k\right)^{1/k} = 1 .
\end{equation*}
\end{itemize}

Then there exists a deterministic positive sequence $(\alpha
_{N})$ subject to $0<\alpha_N \leq N$ for all $N$ and
$\alpha_{N}=o(N)$ as $N\to\infty$, such that
\begin{equation*}
\lim_{N\rightarrow \infty }\frac{1}{N}\nu _{N}\left( \frac{\alpha
_{N}}{N} \right) =1,\quad a.s.
\end{equation*}
and
\begin{equation*}
\lim_{N\rightarrow \infty }\frac{1}{N}\nu _{N}\left( \theta ,\phi
\right) = \frac{\phi -\theta }{2\pi },\quad a.s.
\end{equation*}
In fact the convergence also holds in the $p$\textsuperscript{th} mean for
every positive $p$.
\end{thm}

The proof of this theorem requires the following simple lemma:
\begin{lem}\label{lem:limsup}
If $(\lambda_k)$ is a sequence of real non-negative numbers such
that
\[
\limsup_{k\to\infty} \lambda_k^{1/k} \leq 1
\]
then there exists a sequence $(\varepsilon_N)$ of real positive
numbers tending to zero such that
\[
\sum_{k=0}^N \lambda_k \leq \exp(N \varepsilon_N)
\]
for all $N$.
\end{lem}

\begin{proof}
Let $\epsilon>0$ be arbitrarily small. Since $\limsup_{k\to\infty}
\lambda_k^{1/k} \leq 1$ there exists a constant $C=C(\epsilon)$
such that $\lambda_k \leq C(1+\epsilon)^k$. Therefore
\begin{align*}
\sum_{k=0}^N \lambda_k &\leq C \frac{(1+\epsilon)^{N+1} -
1}{\epsilon}
\\
&\leq \exp\left(N\left(\epsilon + \frac{\log C - \log\epsilon +
\epsilon}{N}\right)\right)
\end{align*}
where we have used the bound $\log(1+\epsilon)\leq \epsilon$ for
$\epsilon>0$. Hence there exists an $N_0=N_0(\epsilon)$ such that
for all $N>N_0$
\[
\sum_{k=0}^N \lambda_k \leq \exp(2N\epsilon),
\]
and since $\epsilon$ is arbitrary, this proves the lemma.
\end{proof}

\begin{proof}[Proof of Theorem~\ref{thm:random polys}]
Note that \eqref{eq:condition on neg power} implies
$\Prob\{|a_{k}|=0\}=0$. Therefore, from
Theorem~\ref{generalclustering} it is sufficient to prove that
there exists a deterministic sequence $\alpha _{N}=o(N)$ such that
$\frac{1}{\alpha_N} L_N \to 0$ a.s., which by Borel-Cantelli would
follow from showing that for any $\epsilon>0$
\[
\sum_{N=1}^\infty \Prob\left\{ \left|\frac{1}{\alpha_N} L_N
\right| \geq \epsilon \right\} < \infty .
\]

First, from the concavity of log, we note that $L_N$ is positive
since
\begin{align*}
L_N &:= \log \left( \sum_{k=0}^{N}|a_{k}|\right) -\frac{1}{2}\log
|a_{0}|- \frac{1}{2}\log |a_{N}| \\
&\geq \log \left( |a_{0}| + |a_{N}|\right) -\frac{1}{2}\log
|a_{0}|- \frac{1}{2}\log |a_{N}| \\
&\geq \log 2 .
\end{align*}
Now, using the fact that
\begin{multline*}
\left\{ \frac{1}{\alpha_N} L_N \geq \epsilon \right\} \\
\subseteq \left\{ \frac{1}{\alpha_N} \log \left(
\sum_{k=0}^{N}|a_{k}|\right) \geq \epsilon/3 \right\} \bigcup
\left\{ \frac{-\log |a_{0}|}{2\alpha_N} \geq \epsilon/3 \right\}
\bigcup \left\{ \frac{-\log |a_N|}{2\alpha_N} \geq \epsilon/3
\right\}
\end{multline*}
we have
\begin{align} \label{eq:show_this_has_finite_sum}
\Prob\left\{ \left|\frac{1}{\alpha_N} L_N \right| \geq \epsilon
\right\} = & \Prob\left\{ \frac{1}{\alpha_N} L_N \geq \epsilon
\right\} \nonumber
\\
\leq &  \Prob\left\{ \frac{1}{\alpha_N} \log \left(
\sum_{k=0}^{N}|a_{k}|\right) \geq \epsilon/3 \right\} 
\\
&+ \Prob\left\{ \frac{-\log |a_{0}|}{2\alpha_N} \geq \epsilon/3
\right\} + \Prob\left\{ \frac{-\log |a_N|}{2\alpha_N} \geq
\epsilon/3 \right\} . \nonumber
\end{align}

To calculate the first sum we wish to bound
\[
\Prob\left\{ \sum_{k=0}^{N}|a_{k}| \geq e^{\epsilon \alpha_N /3}
\right\} .
\]
For this event to happen at least one of the $|a_k|$ must be
bigger than $e^{\epsilon \alpha_N /3 }/ (N+1) $. Therefore,
\begin{align*}
\Prob\left\{ \sum_{k=0}^N |a_k|  \geq e^{\epsilon \alpha_N/3}
\right\}
&\leq \Prob\left\{ \bigcup_{k=0}^N \left\{ |a_k| \geq \frac{1}{N+1} e^{\epsilon \alpha_N /3} \right\} \right\} \\
& \leq \sum_{k=0}^N \Prob\left\{ |a_k| \geq \frac{1}{N+1}
e^{\epsilon \alpha_N/3} \right\} .
\end{align*}
By Tchebychev's inequality we have
\[
\Prob\left\{|a_k| \geq \frac{1}{N+1} e^{\epsilon \alpha_N/3}
\right\} \leq \frac{\E\left[|a_k|^s\right]}{\left(\frac{1}{N+1}
e^{\epsilon \alpha_N/3}\right)^s} = \lambda_k (N+1)^s \exp(-s
\epsilon \alpha_N/3) .
\]
Summing this over $k$ from $0$ to $N$, and using
Lemma~\ref{lem:limsup}, we see that there exists a sequence of
positive numbers $(\varepsilon_N)$ which tends to zero such that
\[
\Prob\left\{ \sum_{k=0}^N |a_k|  \geq e^{\epsilon \alpha_N/3}
\right\} \leq \exp\left(s \log(N+1) - s \epsilon \alpha_N / 3 + N
\varepsilon_N \right) .
\]
Hence if $\alpha_N$ is chosen so that
\[
\alpha_N \geq  \frac{3 N \varepsilon_N}{s \epsilon}+
\left(\frac{3}{\epsilon} +\frac{6}{s\epsilon}\right)\log(N+1)
\]
and subject to the extra conditions $0<\alpha_N \leq N$ with
$\alpha_N = o(N)$, then we see that
\begin{equation}\label{eq:first_sum_bounded}
\sum_{N=1}^\infty \Prob\left\{ \sum_{k=0}^N |a_k|  \geq
e^{\epsilon \alpha_N/3} \right\} \leq \sum_{N=1}^\infty
\frac{1}{(N+1)^2} < \infty .
\end{equation}

We will now deal with the remaining two terms in
\eqref{eq:show_this_has_finite_sum}. If $\alpha_N / \log N \to
\infty$ then by Tchebychev's inequality and the assumption that
$\E\left[|a_0|^{-t}\right] < \infty$, we have
\[
\sum_{N=1}^\infty \Prob\left\{ \frac{-\log |a_{0}|}{2\alpha_N}
\geq \epsilon/3 \right\} \leq \sum_{N=1}^\infty \frac{\E\left[
|a_0|^{-t} \right]}{ \exp(2 t \epsilon \alpha_N / 3)} < \infty .
\]
Finally note that $\limsup_{N\to\infty} \xi_N^{1/N} \leq 1$
implies there exists a positive sequence $(\varepsilon_N')$
tending to zero such that $\xi_N \leq \exp(N \varepsilon_N')$ for
all $N$. Hence by \eqref{eq:condition on neg power} and
Tchebychev's inequality,
\begin{align*}
\sum_{N=1}^\infty \Prob\left\{ \frac{-\log |a_{N}|}{2\alpha_N}
\geq \epsilon/3 \right\} &\leq \sum_{N=1}^\infty \E\left[\frac{1}
{|a_N|^t}\right] \exp(-2 t \epsilon \alpha_N / 3)\\
&= \sum_{N=1}^\infty \xi_N
\exp(-2 t \epsilon \alpha_N / 3) \\
&\leq \sum_{N=1}^\infty \exp\left(N\varepsilon_N'-2 t \epsilon
\alpha_N / 3\right)
\end{align*}
and this sum is finite if $\alpha_N \geq  \frac{3 N
\varepsilon_N'}{2 t \epsilon} + \frac{3}{t\epsilon}\log N$.

Combining the previous two equations and
\eqref{eq:first_sum_bounded} with
\eqref{eq:show_this_has_finite_sum} we conclude that if
$(\alpha_N)$ is a deterministic sequence satisfying $0<\alpha_N
\leq N$ for all $N$, $\alpha_N = o(N)$ and
\[
\frac{\alpha_N }{ \log N + N\varepsilon_N + N\varepsilon_N') } \to
\infty
\]
(and such a sequence exists since both $\varepsilon_N$
and $\varepsilon_N'$ tend to zero), then for any $\epsilon>0$,
\[
\sum_{N=1}^\infty \Prob\left\{ \left|\frac{1}{\alpha_N} L_N
\right| \geq \epsilon \right\} < \infty
\]
which by the Borel-Cantelli lemma shows $L_N / \alpha_N \to 0$
almost surely, and the result now follows from
Theorem~\ref{generalclustering}.
\end{proof}

{}From Theorem~\ref{thm:random polys} we deduce the following
corollary:
\begin{cor}
Let $\left( a_{k}\right) _{k\geq 0}$\ be a sequence of (possibly
dependent) complex random variables such that the moduli $|a_{k}|$
have densities which are uniformly bounded in a neighbourhood of
the origin. Assume that there exists some $s\in \left( 0,1\right]
$ such that if $\lambda_k := \E\left[ \left\vert a_{k}\right\vert
^{s} \right]$ then $\lambda_k < \infty$ for all $k$, and
\begin{equation*}
\underset{k\rightarrow \infty }{\lim \sup }\left( \lambda
_{k}\right) ^{1/k} =1 .
\end{equation*}
Then almost surely the zeros of the classical random polynomial
\begin{equation*}
P(Z) = \sum_{k=0}^N a_k Z^k
\end{equation*}
cluster uniformly around the unit circle.
\end{cor}

\begin{proof}
It suffices to notice that in this special case, $\sup_{N}\xi
_{N}\leq C$ for some positive constant $C$, and so the conclusions
of Theorem~\ref{thm:random polys} follow. These imply uniform
clustering of the zeros (and even give an estimate on the rate of
clustering).
\end{proof}

\begin{example}
Let $P_{N}(Z) =\sum_{k=0}^{N}a_{k}Z^{k}$, with $a_{k}$ being
distributed on ${\mathbb{R}}_{+}$ with Cauchy distribution with
parameter $k^{-\sigma }$, $\sigma >0$. This distribution has
density
\begin{equation*}
\frac{2}{\pi k^{\sigma }}\frac{1}{x^{2}+k^{-2\sigma }}
\end{equation*}
on the positive real line. The conditions of
Theorem~\ref{thm:random polys} are satisfied since $\lambda
_{k}:=\E\left[ a_{k}^{1/2}\right] \leq \frac{C}{k^{\sigma }}$ and
$\xi _{N}:=\E\left[ \frac{1}{a_{N}^{1/2}}\right] \leq Ck^{\sigma
}$. Therefore, if $\alpha _{N}=o(N)$ is such that $\alpha
_{N}/\log N\rightarrow \infty $, then
\begin{equation*}
\lim_{N\rightarrow \infty }\frac{1}{N}\nu _{N}\left( \frac{\alpha
_{N}}{N} \right) =1,\quad a.s.
\end{equation*}
and
\begin{equation*}
\lim_{N\rightarrow \infty }\frac{1}{N}\nu _{N}\left( \theta ,\phi
\right) = \frac{\phi -\theta }{2\pi },\quad a.s.
\end{equation*}
Again, the convergence also holds in the $p$\textsuperscript{th} mean for
every positive $p$.
\end{example}

\section*{Acknowledgments}

This work was partly carried out at the \textit{Random Matrix
Approaches in Number Theory} programme held at the Isaac Newton
Institute, where the first author was supported by EPSRC grant
N09176, partly at the American Institute of Mathematics. Both
authors were partially supported by a NSF Focussed Research Group
grant 0244660. We wish to thank David Farmer, Steve Gonek,
Francesco Mezzadri, Andrew Wade and Marc Yor for useful
conversations. We also wish to thank the referee for careful
reading of this paper and finding a mistake in an earlier draft.

\end{document}